\documentclass[prb,aps,amssymb,footinbib,showpacs]{revtex4-1}
\usepackage{amssymb}
\usepackage{amsmath}
\usepackage{times}
\usepackage{latexsym}
\usepackage{graphicx}
\usepackage{color}
\usepackage{bm}
\usepackage{subfigure}
\usepackage{color}
\usepackage{subfigure}
\usepackage{xcolor}
\usepackage{framed}
\usepackage{mdframed}
\colorlet{shadecolor}{gray!50}

\begin{document}

\title{  Building the Butterfly Fractal: The Eightfold Way}

\author{Indubala I Satija}

\affiliation{Department of Physics and Astronomy, George Mason University, Fairfax, VA 22030}

 \date{\today}

\begin{abstract}

The hierarchical structure of the butterfly fractal -- the Hofstader butterfly,  is found to be described by an octonary tree.  In this framework of building the butterfly graph,  every iteration generates sextuplets of butterflies, each with a tail that is made up of an infinity of butterflies. Identifying {\it butterfly with a tale} as the building block, the tree is constructed with eight generators represented by unimodular matrices with integer coefficients. This Diophantine description provides one to one mapping with the butterfly fractal,
  encoding the magnetic flux interval and the topological quantum numbers of every butterfly.
The butterfly  tree is a generalization of the  ternary tree describing the set of primitive Pythagorean triplets.

 \end{abstract}

\maketitle 

The butterfly graph, with its characteristic ``X" shape structure -  known as the Hofstader butterfly\cite{Hof} is an iconic graph that  describes the energy spectrum of Bloch electrons in a two-dimensional crystalline lattice subjected to a magnetic field. The plot is one of a kind quantum fractal containing some distorted images -- the sub-butterflies ad infinitum.  In this two dimensional landscape, the allowed energies $E$ of the electrons are displayed  as a function of a parameter $\phi$ which specifies the number of flux quanta per unit cell of the lattice. 
The importance of this graph in physics lies in the fact that that it represents a simple model to describe all possible topological states of matter known as the integer quantum Hall states\cite{TKNN}.
The model continues to fascinate physicists as well as mathematicians\cite{ten} since $1960$ when first evidence of the hierarchical nature of spectrum of this system surfaced\cite{azbel,others}.
 Furthermore, there have been various attempts to observe the intricate spectral structure in a real laboratory setup\cite{Dean}.\\

Butterfly graph is a goldmine for number theories as it is perhaps the only example in physics literature where number theory finds application\cite{RF}.  Recent studies have dwelled on its geometrical and number theoretical aspects\cite{book,SAT16, SW, SAT21, SAT21E} revealing its relationship to other abstract hierarchical structures such as the Farey tree, the Apollonian gaskets, the tree of primitive Pythagorean triplets. Most of these studies have focussed on characterization of self-similar features of the fractal.\\

In this paper, we present an iterative scheme to build this fractal, visualized as an octonary  tree. This ``Butterfly tree" has  eight generators, where six of them generate the six butterflies and the other two generate the chain ( left and right chain ), consisting of infinity of butterflies that is attached to the butterfly like a tale. In its simplest representation, each generator  of the group, which we refer as the ``Butterfly group",  is a 
 $3\times3$ matrix with integer entries, acting on a column matrix consisting of $3$ integers that uniquely label a butterfly. This representation encodes its magnetic flux interval and the topological quantum numbers of the two major gaps that label every butterfly. This eightfold way of generating the butterfly graph where the unique labelling of the butterfly and their recursions are  described  in terms of integers is a  ``Diophantine representation" of the butterfly problem.\\
 
 Section I summarizes some  previously known results such as the Wannier diagram\cite{W}, the  the Farey and their relation to the butterfly graph. It  shows how it leads to  three integer that
 uniquely label every  butterfly in the butterfly graph. Section II discusses the eightfold way of building the butterfly fractal. Here we first discuss the tree of Pythagorean triplets that describes the sub-structure of the butterfly tree characterized as parity conserving as well as conserving the asymmetry of the butterfly\cite{PT1}. We then discuss the parity violating part of the tree which does not preserve symmetry characterization.  We determine the eight generators that provide an iterative scheme to build the butterfly fractal.
 Section III  discusses the relationship between the Butterfly group and the super Apollonian group.
 \section{ Planting the Butterfly Tree}

To set the stage for describing  the butterfly tree, we present a brief review of some of the  previously discussed number theoretical aspects of the butterfly fractal, along with some new perspective
to this problem.

\subsection{ Wannier Diagram,  Farey Tree and their Relationship}

Butterfly graph has a simple representation known as the Wannier diagram\cite{W} that provides a simple representation of the spectrum by labeling all the gaps of the spectrum with two integers $(\sigma, \tau)$  expressed as a linear Diophantine equation\cite{Dana},
\begin{equation}
 p \sigma + q \tau = r
\label{Dio}
\end{equation}
  Here $r$ labels the $r^{th}$ gap of the spectrum for a rational magnetic flux $\phi = \frac{p}{q}$. Here $\sigma$ is the Chern number - the topological quantum number of the gap. The Chern number $N$ of the band also satisfies a Diophantine equation,
  \begin{equation}
 p N + q M = 1
\label{Diob}
\end{equation}

\begin{figure}[htbp] 
 \includegraphics[width = .33 \linewidth,height=.42 \linewidth]{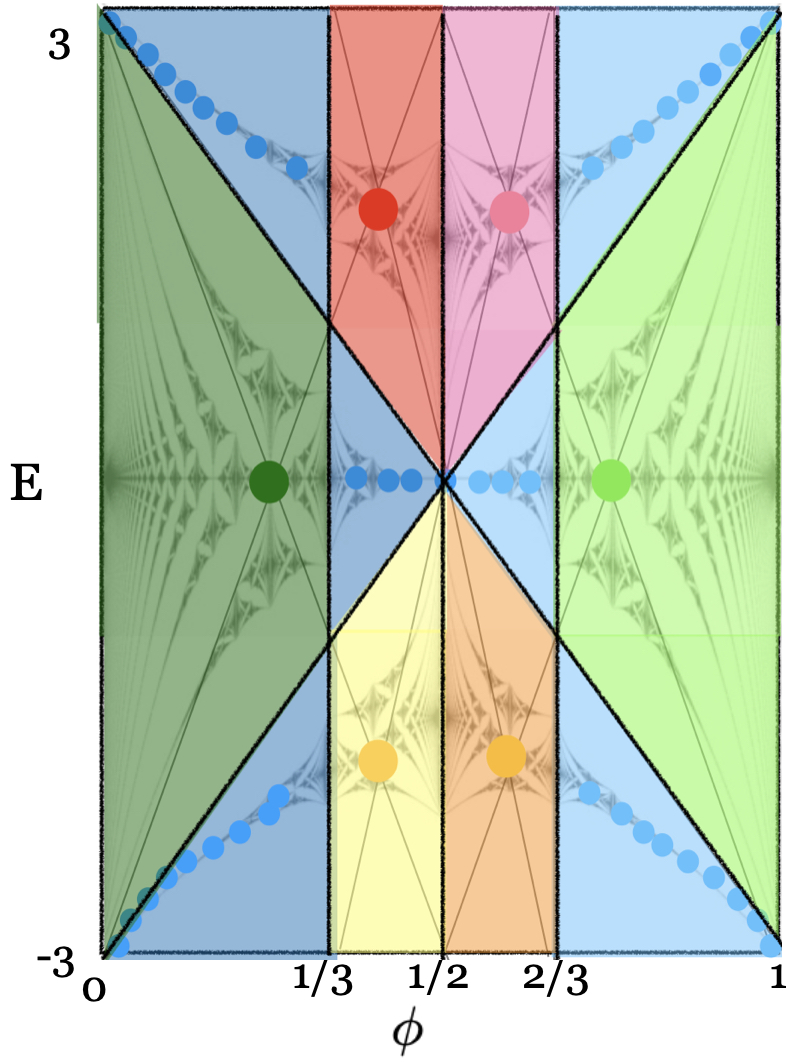} 
\caption{  Figure shows the skeleton butterfly graph overlayed on the butterfly fractal. The six trapezoids in the graph house the six butterflies represented by a dot at the center of the butterfly.
Attached to one of the  parallel lines of every trapezoid is a triangle ( shown in blue ) that overlays an infinite chain of butterflies. The eight shades of colors used in this graph correspond to eight generators needed to produce this sub-image from the main butterfly as will be explained in next section.}
\label{ball}
\end{figure}

The Wannier diagram can be viewed as representing  ``butterfly skeleton" as the channels of forbidden energies shrink to straight lines. The skeleton graph is a geometrical construction,  a kind of tiling of trapezoids and triangles where all slanting lines represent the gaps of the butterfly fractal.
Fig. (\ref{ball}) displays color coded skeleton graph overlayed on the butterfly fractal.
 This figure is a reorganization of the Wannier diagram with six trapezoids, representing six butterflies, each attached to a triangle consisting of a chain of monotonically
 decreasing sizes of butterflies that degenerate to a point either to the left or to the right of the butterfly. Rationale for this way of assembling the skeleton butterfly and labeling it with eight shades of colors is a signal to the eightfold way of constructing the butterfly tree, as will be discussed shortly.\\ 

The key feature of the two-dimensional landscaping in the skeleton graph with trapezoids and triangles is the fact that the diagonals of the trapezoid as well as all the slanting lines in the graph, all representing gaps in the butterfly graph,
 have integer slopes. This feature is rooted in the fact that the two parallel lines of all the trapezoid are Farey neighbors\cite{SAT21}. That is, the vertical lines at the rational fractions $\frac{p_L}{q_L}$ and $\frac{p_R}{q_R}$  satisfy the following relation:

\begin{equation}
p_L q_R - p_R q_L = \pm 1,
\label{Flr}
\end{equation}

where  $(p_L ,q_L)$ and $(p_R ,q_R)$ are coprimes. Any two fractions satisfying  Eq. (\ref{Flr}) are two neighboring fractions in the Farey tree and are known as the {\it friendly fractions}.  Farey tree is constructed by applying the `` Farey sum rule" to  $\frac{p_L}{q_L}$ and $\frac{p_R}{q_R}$- the Farey parents that gives a new  fraction $\frac{p_c}{q_c}$ - the Farey child:
\begin{equation}
\frac{p_c}{q_c} = \frac{p_L}{q_L} \bigoplus \frac{p_R}{q_R} =  \frac{p_L+q_R}{q_L+q_R}.
\label{Fc}
\end{equation}
Analogous to the friendly pair $\frac{p_L}{q_L}$ and $\frac{p_R}{q_R}$, $\frac{p_c}{q_c}$ also forms a friendly pair with each of its parents $\frac{p_L}{q_L}$ and $\frac{p_R}{q_R}$.
These equations  define a Farey triplet denoted as $\Big[\frac{p_L}{q_L}, \frac{p_c}{q_c}, \frac{p_R}{q_R}\Big]$ which will be referred  as the ``{\it friendly Farey triplet}". The $x$-coordinate of the point of intersection of the trapezoid  coincides precisely with the $x$-coordinate of the center of the butterfly in the butterfly graph.\\

In the butterfly graph, the three fractions in the 
Farey triplet have been shown to describe the left boundary, the center and the right boundary of a butterfly\cite{book,SW, MW1,MW2}, showing an intimate link between the number theory and the quantum problem of competing lengths which lies at the heart of the butterfly fractal. The topological quantum numbers of the butterfly emerge simply from number theory as they are
the integer slopes of the diagonals of the trapezoids.  For two major gaps of a butterfly, the corresponding two Chern numbers will be denoted as
$(\sigma_+, - \sigma)$ where $(\sigma_+, \sigma_-)$ are both positive.
With the exception of the trapezoids in the central part of the energy spectrum, symmetric about $E=0$,  the other trapezoids in the skeleton graph do not exhibit horizontal mirror symmetry. As described in our recent study\cite{SAT21E}, all  the trapezoidal cells representing butterflies are governed by the rule of 
 {\it minimum violation of horizontal mirror symmetry}.  This rule is quantified by the parameter $\Delta \sigma$ defined as,
 
 \begin{equation}
 \Delta \sigma = \sigma_+ -  \sigma_- .
 \end{equation}
 
In other words, from  the $q_L$ bands located at $\frac{p_L}{q_L}$ and $q_R$ bands located at $\frac{p_R}{q_R}$, one of the bands from each set mates with the other set to form a butterfly provided they abide by the rule of minimal violation of symmetry, in addition to being Farey neighbors. As shown below, eightfold way of constructing butterfly graph has both these rules embedded in the eight generators that produce the butterfly graph.
  
  \subsection{ Dichotomy in the butterfly graph - The C-cell and the E-cell butterflies }

 The two main gaps, forming the X shaped structure divide the butterfly into central and edge parts, referred as the C-cell and the E-cells\cite{SAT21E} of the butterfly. This geometrical division actually has a much deeper significance. Firstly, C-cell butterflies conserve parity while the E-cell butterflies do not. The parity, described by the parity of $q_c$, the denominator of the rational flux value corresponding to the center of the butterfly. Furthermore, this dichotomy is also reflected in the topological characterization of a  butterfly as the  C-cell  butterflies  preserve $ \Delta \sigma$, the asymmetry parameter defined above while $\Delta \sigma$ grows exponentially in the E-cell butterflies. 
Furthermore, although not obvious from the geometrical division,  the chain of butterflies attached to every  butterfly are part of the C-cell butterflies as all members of the chain have same parity
and same $\Delta \sigma$. These two distinctive features of the butterfly and the important role they play in constructing the butterfly tree will be discussed in detail in next section.

  \subsection{ Butterfly Integers: Diophantine Representation of the Butterfly plot}

In the butterfly graph, six positive integers $( p_L, p_R, q_L, q_R, \sigma_+, \sigma_-)$ are associated with a butterfly, providing a complete description of the skeleton butterfly represented by a trapezoid in the skeleton diagram. This mapping is an isomorphism between the skeleton and the actual butterfly in the fractal graph, in sharp contrast to the homomorphism  mapping 
described earlier\cite{SW} by $( p_L, p_R, q_L, q_R, M , N)$ where $(M,N)$  describe the band-topology  of the butterfly. \\

We now argue that just three integers -  the triplet $(q_R, q_L, \Delta \sigma)$  are sufficient to label a butterfly uniquely as we impose the Farey rules given by Eqs (\ref{Flr}) and (\ref{Fc}). These three integers  determine the magnetic flux and topological quantum numbers  of a butterfly. Firstly, the  $(q_R, q_L)$ determine the magnetic flux values $\Big[ \frac{p_L}{q_L}, \frac{p_c}{q_c}, \frac{p_R}{q_R} \Big ]$ in view of the Diophantine equation. That is, 
 given a pair $(q_R, q_L)$, the corresponding  pair $(p_R, p_L)$ is uniquely determined from the Diophantine equation $ (p_L q_R - p_R q_L ) = \pm 1$ as both $(p_L, q_L)$ and $(p_R,q_R)$ are relatively prime and the fractions $\frac{p_L}{q_L}$ and $\frac{p_R}{q_R}$ are less than unity. Secondly, given $\Delta \sigma$, we know $(\sigma_+, \sigma_-)$ as
 $\sigma_+ + \sigma_- = q_c = q_R+q_L$. It should be noted that the pair $(q_R,q_L)$ has handedness built into it as $\frac{p_R}{q_R} ( \frac{p_L}{q_L}) $ are always taken as the right (left) boundary of the butterfly. \\
 
\section{ Eight-fold-Way}

  \begin{figure}[htbp] 
 \includegraphics[width = .65 \linewidth,height=.35 \linewidth]{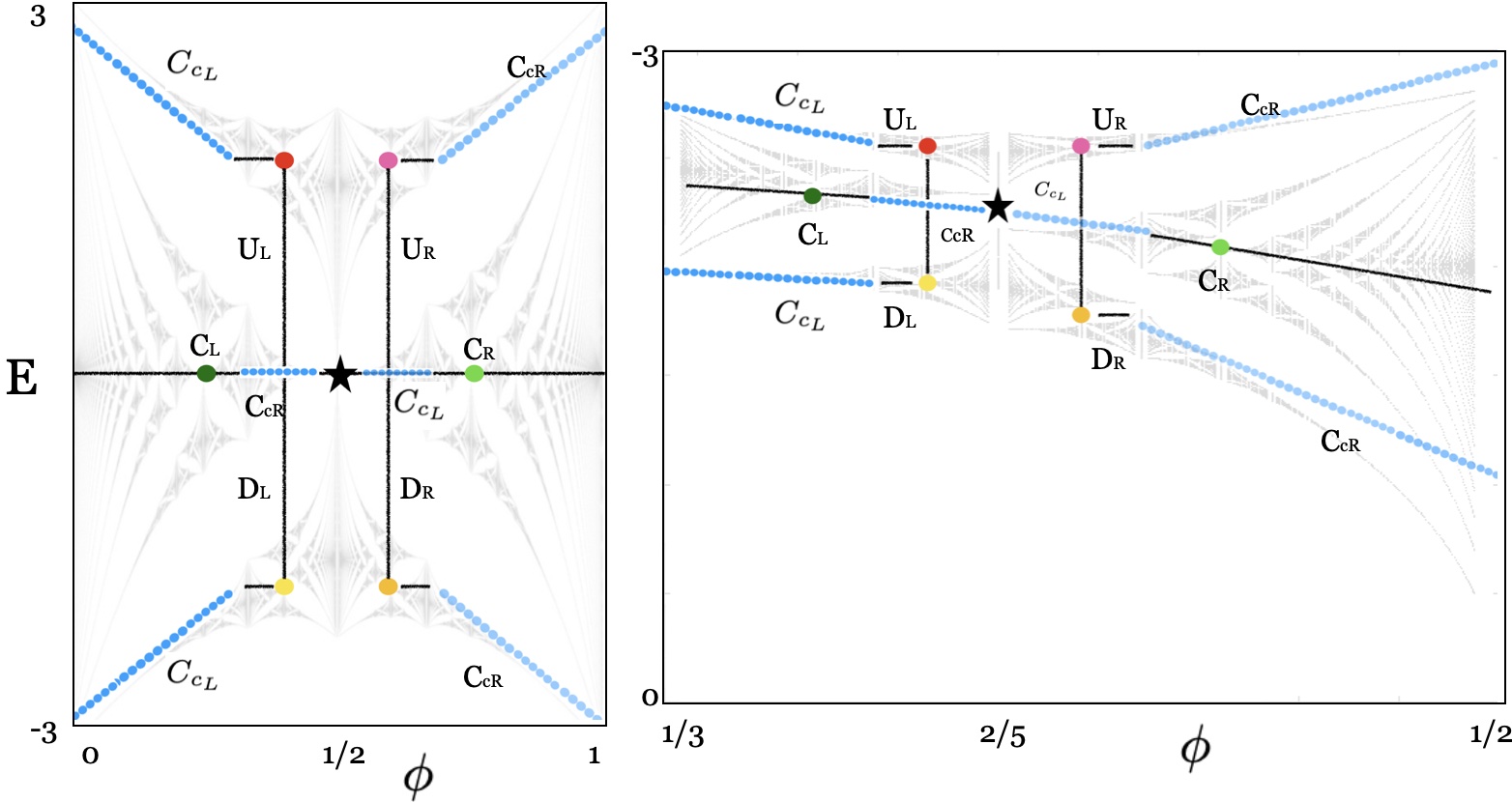} 
\caption{ Eightfold partitioning of the butterflies, shown with eight colors where six butterflies are shown with a  dot at their center and left and right chains (shown respectively with two shades of blue ) displayed with a line of dots. Each color is associated with a generator labelled explicitly. The chains consist of infinity of butterflies with monotonically decreasing sizes,  terminating at a point which is either to the left or to the right of the butterfly from which the chain originates. This associates a kind of handedness to the description of the chains, requiring two distinct generators, one for the left and the  other for the right chain. Left and right panels in the figure show the main butterfly and one of its sub-images. }
\label{b8}
\end{figure}

We will now describe the eightfold way of constructing the butterfly fractal. Fig. (\ref{b8}) summarizes the key aspects of the theoretical framework.
We  label the sub-structure with eight generators.
Six of these generators produce six baby butterflies. Each infant is  attached to an entourage- a chain consisting of infinite number of butterflies. In this dissection, three infants reside to the left and the other two to the right of the center of the parent  butterfly. The generators for the two C-cell butterflies are denoted as $C_L$ and $C_R$. The four generators of the E-cell butterflies are labelled as $(U_L, U_R)$  ( for upper butterflies ) and $(D_L, D_R)$ for the lower butterflies. Generators for the right and the left chains are labelled as $C_{cR}$ and $C_{cL}$ as the butterflies in the chains preserve parity and also $\Delta \sigma$ - the two key features that define the C-cells.
Note that although chains consist of infinity of butterflies, all the members of a chain are generated  recursively by a single generator, and thus justifying the same color to describe all the butterflies in the chain.The resulting octet tree describes the iterative process of constructing the skeleton butterfly graph.\\

We will first focus on the central part of the tree, describing butterflies that exist in the C-cells. Earlier studies have suggested that this subset of the recursive structure which is well described by the tree of primitive Pythagorean triplets\cite{SAT21,PT1}. Below we review this result, highlighting some subtle differences between the Pythagorean tree and the tree describing the C-cell butterflies.

 \subsection{ Pythagorean Tree:  Parity Conserving branches of the Butterfly Tree }
 
 \begin{figure}[htbp] 
   \includegraphics[width = .5 \linewidth,height=.4 \linewidth]{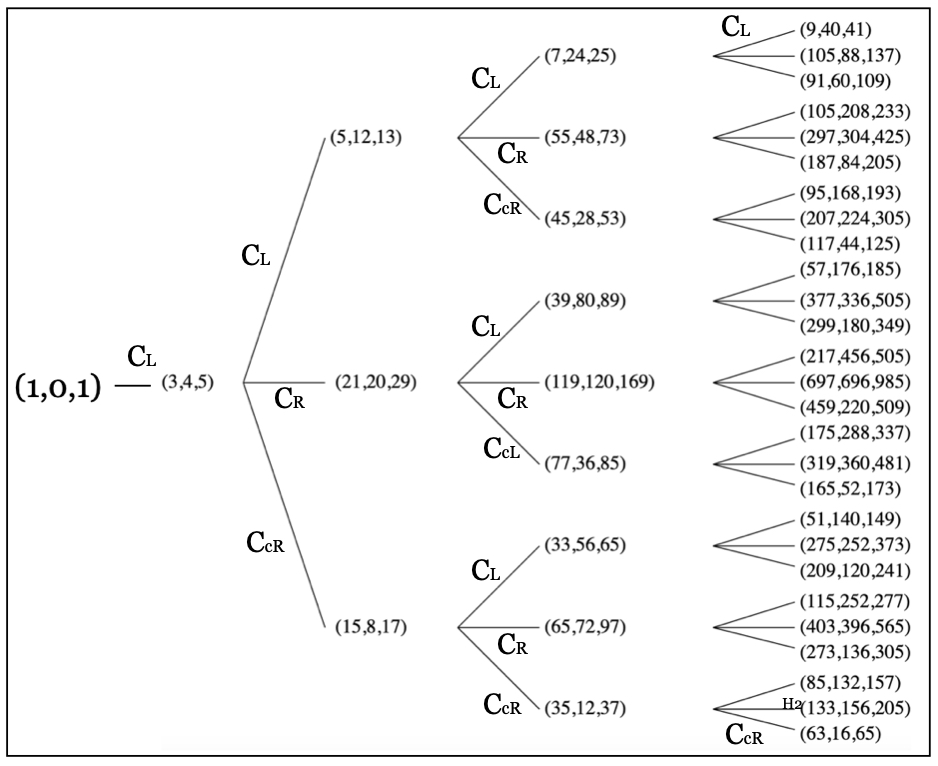} \quad
 \includegraphics[width = .7 \linewidth,height=.4 \linewidth]{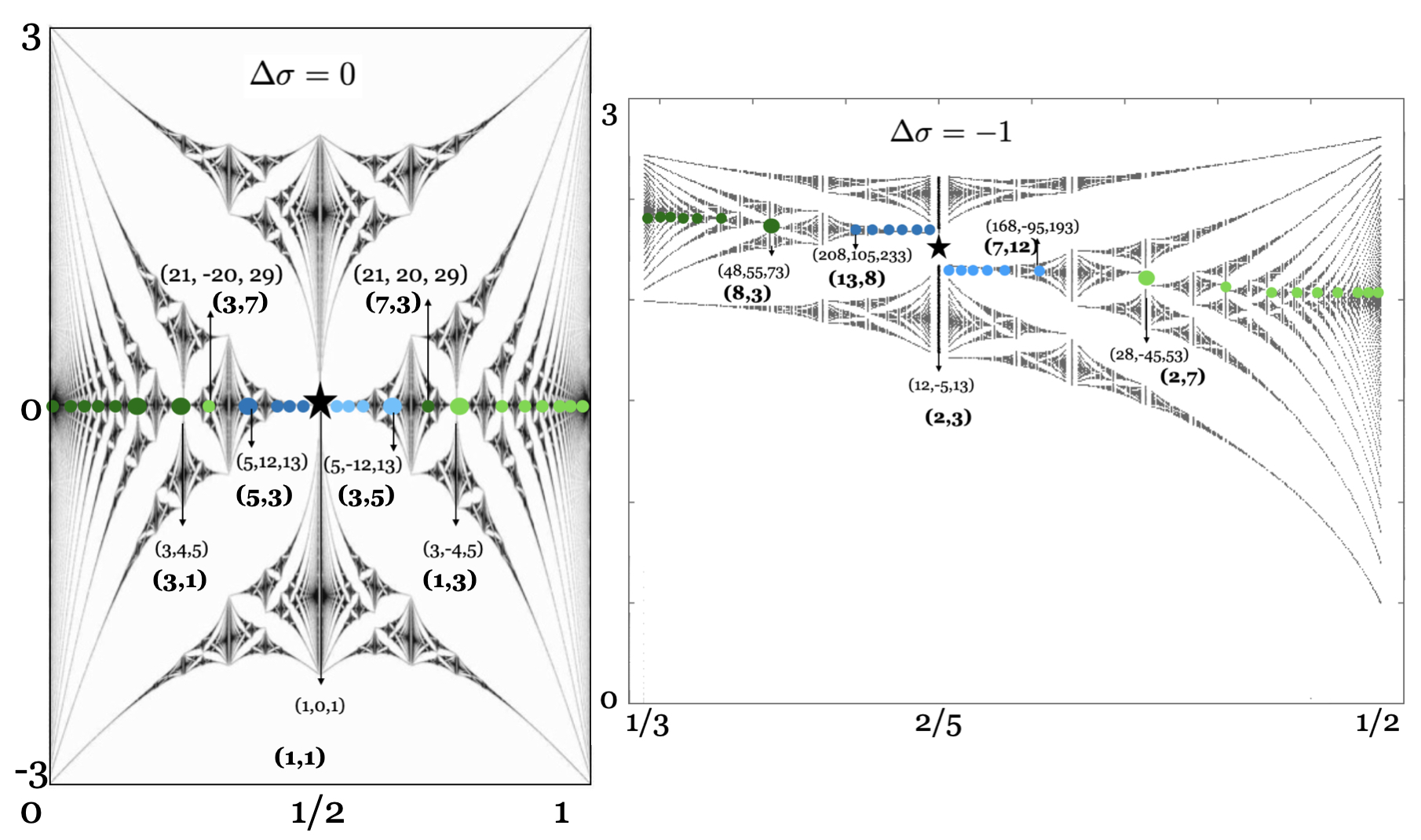} 
 \caption{  Upper panel shows the Pythagorean tree - an abstract structure labelled with  three generators. Lower panels show the main butterfly ( left )  and one of its sub-image ( right ), with their distinct $\Delta \sigma$ labels. The sub-butterflies of the C-cells of these two butterflies are color coded  where dark and light green dots respectively show the $H_1$ and $H_2$ hierarchies and the blue shows the corresponding chains generated by $H_3$.  Some of the sub-images are labelled with Euclid parameters ( in bold ) and Pythagorean triplets so one can identify the Pythagorean tree hierarchy in  these butterflies.}
\label{PT}
\end{figure}
 
A Pythagorean triple is a set of three positive integers $(a, b, c)$  satisfying the equation,

\begin{equation}
 a^2+b^2=c^2.
 \label{P}
 \end{equation}
 
The triplet is said to be primitive if and only if a, b, and c share no common divisor.  
More than two thousand years ago, Euclid provided a recipe for generating such a triplets of integers with
 a pair of integers $(m,n)$ which we will refer as  the {\it Euclid parameters}. It generates a Pythagorean triplet, given by,

\begin{eqnarray*}
(a,b,c)  & =   & \Big( mn ,\,\ \frac{m^2-n^2}{2}, \frac{m^2+n^2}{2} \Big)  \,\,\ \rm{ if \,\  m, \,\ n \,\ have \,\  same \,\ parity }\\
& =   & \Big ( 2mn , m^2-n^2,  m^2+n^2 \Big)   \,\,\ \rm{ if \,\  m, \,\ n \,\ have \,\  opposite \,\ parity } .
\end{eqnarray*} 

It has been known that \cite{BPT}, \cite{Hall}  the entire set of all primitive Pythagorean triples has the structure of a tree. 
The tree  as shown in (\ref{PT})  is generated iteratively by three matrices $H_1$, $H_2$, $H_3$, acting on the column matrix with entries as the Pythagorean triplet $(a,b,c)$:
\begin{equation}
H_1= \left( \begin{array}{ccc} 1 & -2 & 2  \\    2 & -1 & 2  \\  2 & -2  & 3 \\ \end{array}\right),
H_2=\left( \begin{array}{ccc} 1 & 2 & 2  \\    2 & 1 & 2  \\  2 & 2  & 3 \\ \end{array}\right),
H_3=\left( \begin{array}{ccc} -1 & 2 & 2  \\    -2 & 1 & 2  \\  -2 & 2  & 3 \\ \end{array}\right) 
\label{H3}
\end{equation}

We can simplify the characterization of the tree  using
 $ 2\times 2$ matrices, acting on the pair of Euclid parameters $(m,n)$,  

\begin{equation}
h_1 =  \left( \begin{array}{cc} 1   & 2  \\  0 & 1   \\  \end{array}\right), \quad  h_2 =  \left( \begin{array}{cc} 2   & 1  \\  1 & 0   \\  \end{array}\right),\,\,\,\ h_3 =  \left( \begin{array}{cc} 2   & -1  \\  1 & 0   \\  \end{array}\right)
\end{equation}

Hidden in these matrices is one of the key feature of the tree, namely the parity. It can be quantified by  $q_c=q_R+q_L$. Even $q_c$ corresponds to both $(m,n)$ being odd and odd-$q_c$ corresponds to $m$ and $n$ having opposite parity. This  can be see from the following equations:

\begin{eqnarray}
\rm { with} \,\ h_1 : q'_c  &  =   &  q_c+ 2 q_{L} ,\,\ \rm { : preserves \,\ q_c \,\ parity } \nonumber \\
\label{h123}
\rm { with} \,\ h_2 : q'_c    & = &  q_c+2 q_{ R}  \,\ \rm { : preserves \,\ q_c  \,\ parity } \nonumber \\ 
\rm { with} \,\  h_3  : q'_c    & = &  -q_c + 4 q_{R}  \rm { : preserves  \,\  q_c \,\ parity } \nonumber\nonumber
\end{eqnarray}

Parity conservation is also reflected in the Pythagorean triplets as the parity ( even or odd feature ) of the each member of the triplet is conserved under the application of $H_i$.\\

It turns out that if we use butterfly integers $(q_R, q_L)$ or $( p_R, p_L)$ as Euclid parameters, the Pythagorean tree almost mimics all the C-cell hierarchies of the butterfly graph.
The reason why only C-cell recursions follow the Pythagorean tree pattern is due to the fact that only the C-cell butterflies conserve parity ( see below). As a note of caution, we cannot use $(p_x, q_x)$
( x = L, R ) as Euclid parameters to describe butterfly hierarchy with the Pythagorean tree as a given pair $( p_x, q_x)$ does not determine a unique C-cell butterfly.\\

 To relate the butterfly recursions with the Pythagorean tree, we need two additional features to be included in the tree structure.\\
   
(I)  \underline{ Allowing negative entries in the Pythagorean triplets }\\

 The Pythagorean tree is generated using Euclid parameters $(m,n)$ where $ m > n$. This results in three Pythagorean triplets to be positive integers.
For the C-cell butterflies, with $(q_R,q_L)$  as the Euclid parameters, there are butterflies where $q_R > q_L $ as well as those where $q_R < q_L$.
This results in triplets with  $ b > $ or $b < 0$. In other words, Pythagorean triplet representing a C-cell butterflies will also have negative integers in the set as shown in the left panels in Fig. (\ref{PT}) \\

To accommodate these two possibles, we need to introduce four matrices. It is better to work with $2\times2$ matrices that generate new set of Euclid parameters and hence relate more directly to the butterfly labeling. Relabeling the matrices to represent the left and the right  babies of the parent butterfly, the four matrices are:  

\begin{equation}
c_L =  \left( \begin{array}{cc} 1   & 2  \\  0 & 1   \\  \end{array}\right), \quad c_R =  \left( \begin{array}{cc} 1   & 0  \\  2 & 1   \\  \end{array}\right),  \,\  \quad  c_{cR}=  \left( \begin{array}{cc} 2   & -1  \\  1 & 0   \\  \end{array}\right) , \,\ \quad c_{cL} =c^{-1}_{c R} =\left( \begin{array}{cc} 0   & 1  \\  -1 & 2   \\  \end{array}\right)
\end{equation}

We note that $c_{cL}$ and $c_{cR}$ generates butterfly chains as 
${\rm with \,\ c_{cL} } : q_L \rightarrow q_R$ and 
${\rm with \,\ c_{cR} } : q_R \rightarrow q_L$. Both these operators conserve parity.\\

(II) \underline{Accommodating butterfly siblings}\\

\begin{figure}[htbp] 
 \includegraphics[width = .48 \linewidth,height=.55 \linewidth]{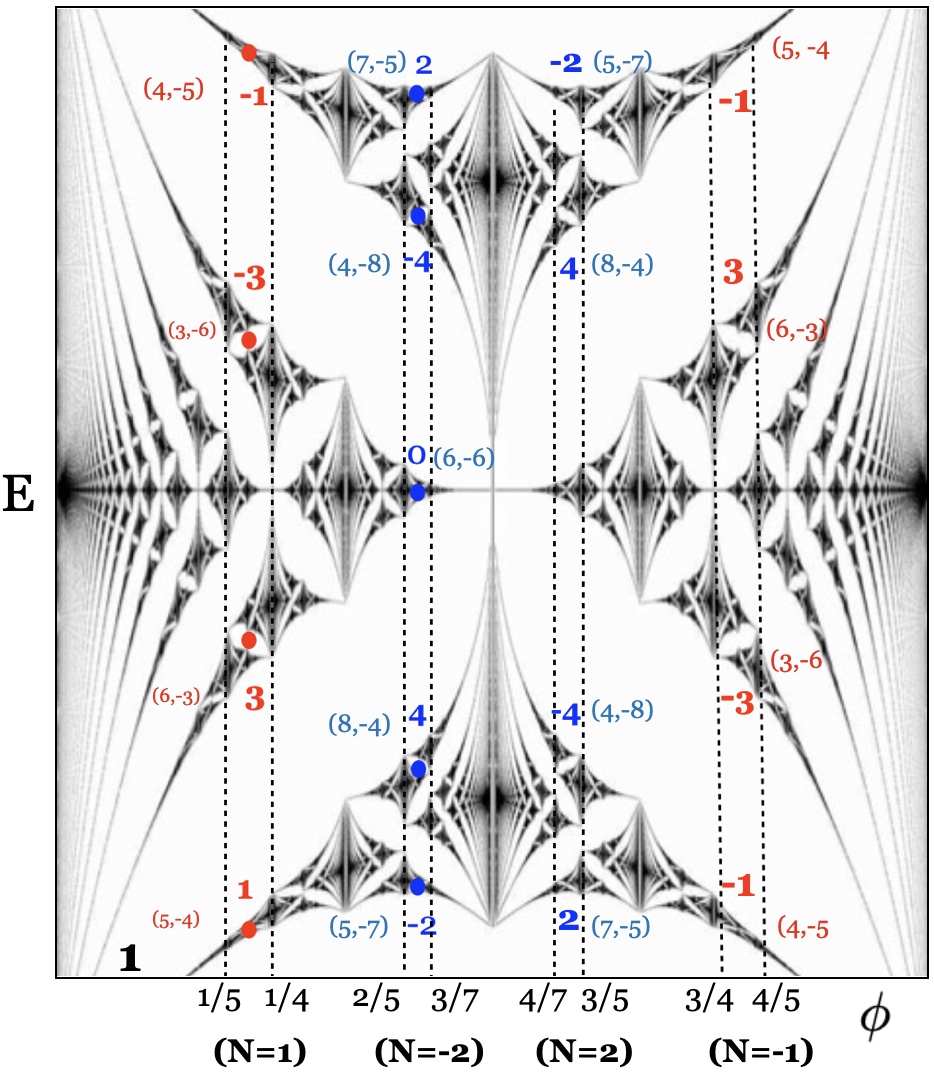} \quad
\caption{ Four families of butterflies, showing butterfly siblings,  the butterflies  that share the same flux interval . Within a family, all the siblings are labeled with same Euclid parameters and also the band Chern number $N$ ( shown in bold black). They are are distinguished with $\Delta \sigma$ ( in bold)-  color coded  integers  and $(\sigma_+,\sigma_-)$ shown in the figure. Note that ``dots" labeling butterfly centers are omitted from the $\phi > 1/2$ part of the to avoid obstruction in displaying the fine structure of the fractal.}
\label{sib}
\end{figure}

Family of butterflies that reside in the same flux interval, appear vertically stacked in the butterfly graph have their members that are topologically distinct. Earlier study\cite{SW} where the self-similar recursive pattern of the butterflies was described by integers $(p_L, p_R, q_L, q_R, M, N)$ does not distinguish the siblings as they not only share the magnetic flux interval but also the topological integers $(M,N)$ where $N$ is the Chern number of the  boundary bands of the butterflies. Fig. (\ref{sib}) shows explicitly that siblings are distinguished by $(\sigma_+, \sigma_-)$ or by $\Delta \sigma$.  This clearly shows that it is the quantum numbers $\sigma_{\pm}$  that capture the topology of the butterfly and not the band Chern number $N$. 

For the central butterflies, that is butterflies with horizontal mirror symmetry, $\Delta \sigma=0$ as $\sigma_+ = \sigma_- = \frac{q_c}{2}$. For other C-cell butterflies,
this mirror symmetry is lost. The recursion relations for the topological integers $(\sigma_+, \sigma_-, \Delta \sigma)$ for these butterflies  will be described in the next sub-section which proves that
for C-cell butterflies, $\Delta \sigma$ remains a constant under renormalization.\\

We close this section by emphasizing that the infinity of  butterflies in a given chain ( shown in blue ) are part of the C-cell of the butterfly to which the chain is attached.
Each chain is like a ``tail" attached a butterfly . A Chain carries over the parity and also the
$\Delta \sigma$  of the butterfly that it is attached to. All members of a given chain are generated by a the same generator and this generator and its inverse produce the left and the right chain.
This  is the rationale for partitioning the skeleton graph into six butterflies with a tail. 

\subsection{  Parity  violating branches of the butterfly tree}

A closer look at the butterfly fractal reveals recursive pattern for the for E-cell butterflies as seen in Figs (\ref{ball},\ref{b8}). Given  a parent butterfly $\Big[\frac{p_L}{q_L}, \frac{p_c}{q_c}, \frac{p_R}{q_R} \Big ]$,   a simple observation leads to the following recursions.\\

(I) The butterfly  on the left( upper and lower ) is obtained from the parent butterfly using
\begin{equation}
  \frac{p_c}{q_c} \rightarrow  \frac{p_R}{q_R}  \,\,\ \frac{p_L}{q_L} \rightarrow  \frac{p_c}{q_c} \bigoplus \frac{p_L}{q_L} = \frac{p_c+p_L}{q_c+q_L}
  \label{l1}
  \end{equation}

(II) The right  butterfly (upper and lower) is obtained from the parent butterfly using ,
\begin{equation}
 \frac{p_c}{q_c} \rightarrow  \frac{p_L}{q_L}  \,\,\ \frac{p_R}{q_R} \rightarrow  \frac{p_c}{q_c} \bigoplus \frac{p_R}{q_R} = \frac{p_c+p_R}{q_c+q_R}
 \label{l2}
 \end{equation}

Translated in the matrix form, we have two new $2\times 2$ matrices that act on   $(q_R, q_L)$ ( or $(p_R, p_L)$ ) to produce new sets of these butterfly integers:
\begin{equation}
u_L = d_L =  \left( \begin{array}{cc} 1   & 1  \\  1 & 2   \\  \end{array}\right), \quad  u_R = d_R =  \left( \begin{array}{cc} 2   & 1  \\  1 & 1   \\  \end{array}\right),
\end{equation}
We can see that these matrices lead to hierarchies that do not conserve parity:

\begin{eqnarray}
\rm { Using} \,\  u_L: q'_c    & = &  2q_c+ q_{ L}  \,\ \rm { : does \,\ not  \,\ preserves \,\ q_c  \,\ parity } \nonumber \\ 
\rm { Using} \,\ u_R : q'_c    & = &  2q_c + q_{R}  \rm { : does \,\  not \,\ preserves  \,\  q_c \,\ parity } \nonumber\nonumber
\end{eqnarray}

Consequently, these hierarchies cannot be described by the Pythagorean tree. In addition, we also need the recursive scheme for the topological quantum numbers
$(\sigma_+, \sigma_-)$ or $\Delta \sigma$. \\

To obtain the recursion relations for $(\sigma_+, \sigma_-)$, we make use of the result\cite{book} that the Chern numbers of the hierarchical set of gaps near a rational flux value say $\frac{p_0}{q_0}$ are $\sigma_0+ nq_0$ where $\sigma_0$ is the Chern number of the parent gap.  As an example to illustrate this,  let us determine the Chern numbers of the major gaps of the edge butterfly residing 
in the interval $\frac{1}{3}, \frac{2}{5}, \frac{1}{2}$ ( shown in Figs (\ref{ball}, (\ref{b8}) ) with parent butterfly being the main butterfly $\frac{0}{1}, \frac{1}{2}, \frac{1}{1}$ with Chern numbers $(1,1)$. 
The gap near $\frac{p_0}{q_0} = 1/3$ has Chern number $\sigma_0=-1$. Therefore, the Chern numbers $(\sigma_+, \sigma_-)$ of the infant butterfly are $ -1+3=2$ and $-1-2=-3$.\\

Using the equations (\ref{l1}) and (\ref{l2}), along with  the identity, $\sigma_+ + \sigma_- = q_c$, we can derive the renormalization trajectories of the
Chern numbers $( \sigma_+, \sigma_-)$ . The results are summarized in Table 1 and are also given in eight matrices below. These matrices  acting on the column matrix consisting of four butterfly integers $(  q_R, q_L, \sigma_+, \sigma_-)$ are:

 \begin{eqnarray*}
\mathcal{U}_L  & =  & \left( \begin{array}{cccc}    1 & 1 & 0 & 0 \\   1 & 2 & 0 & 0 \\    0 & 1 & 1 & 0 \\   1 & 1 & 0 & 1 \\
\end{array}\right),\,\,\ 
\mathcal{U}_R = \left( \begin{array}{cccc}    2 & 1 & 0 & 0 \\  1 & 1 & 0 & 0 \\    1 & 1 & 1 & 0 \\   1 & 0 & 0 & 1 \\
\end{array}\right), \,\,\,\ 
\mathcal{D}_L = \left( \begin{array}{cccc}  1 & 1 & 0 & 0 \\   1 & 2 & 0 & 0 \\    1 & 1 & 1 & 0 \\   0 & 1 & 0 & 1 \\
\end{array}\right),\,\,\ 
\mathcal{D}_R = \left( \begin{array}{cccc}   2 & 1 & 0 & 0 \\  1 & 1 & 0 & 0 \\   1 & 0 & 1 & 0 \\    1 & 1 & 0 & 1 \\
\end{array}\right) \\
\\
\mathcal{C}_L  & =  & \left( \begin{array}{cccc}   1 & 2 & 0 & 0 \\   0 & 1 & 0 & 0 \\    0 & 1 & 1 & 0 \\   0 & 1 & 0 & 1 \\
\end{array}\right),\,\,\ 
\mathcal{C}_R = \left( \begin{array}{cccc}    1 & 0 & 0 & 0 \\   2 & 1 & 0 & 0 \\    1 & 0 & 1 & 0 \\   1 & 0 & 0 & 1 \\
\end{array}\right), \,\,\,\ 
\mathcal{C}_{cL} = \left( \begin{array}{cccc} 0 & 1 & 0 & 0 \\   -1 & 2 & 0 & 0 \\    1 & -1 & 1 & 0 \\   1 & -1 & 0 & 1 \\
\end{array}\right),\,\,\ 
\mathcal{C}_{cR} = \left( \begin{array}{cccc}   2 & -1 & 0 & 0 \\  1 & 0 & 0 & 0 \\    1 & -1 & 1 & 0 \\   1 & -1 & 0 & 1 \\
\end{array}\right)
\end{eqnarray*}

\bigskip\par\noindent

As discussed earlier,  just three integers are sufficient to label a butterfly, which we choose to be $(q_R, q_L, \Delta \sigma$).  Therefore, the simplest representation of the eight generators are $3 \times 3 $ matrices which are given by:
\begin{eqnarray*}
U_L & = & \left( \begin{array}{ccc} 1 & 1 & 0  \\   1 & 2 & 0 \\  0 & 1 & 1  \\ \end{array}\right), \,\,\,\ U_R= \left( \begin{array}{ccc} 2 & 1 & 0  \\   1 & 1 & 0 \\  0 & 1 & 1  \\ \end{array}\right),\,\,\ 
D_L  = \left( \begin{array}{ccc} 1 & 1 & 0  \\   1 & 2 & 0 \\  0 & 1 & 1  \\ \end{array}\right) ,\,\,\ D_R = \left( \begin{array}{ccc} 2 & 1 & 0  \\   1 & 1 & 0 \\  0 & 1 & 1  \\ \end{array}\right)\\
\\
C_L  & = & \left( \begin{array}{ccc} 1 & 2 & 0  \\   0 & 1 & 0 \\  0 & 0 & 1  \\ \end{array}\right) ,\,\,\ C_R =\left( \begin{array}{ccc} 1 & 0 & 0  \\   2 & 1 & 0 \\  0 & 0 & 1  \\ \end{array}\right)  ,\,\,\ C_{cL} = \left( \begin{array}{ccc} 0 & 1 & 0  \\   -1 & 2 & 0 \\  0 & 0 & 1  \\ \end{array}\right) ,\,\,\ C_{cR} =\left( \begin{array}{ccc} 2 & -1 & 0  \\   1 & 0 & 0 \\  0 & 0 & 1  \\ \end{array}\right) 
\end{eqnarray*}

Table 1 lists the recursion relations and shows examples that illustrate the scheme for constructing the tree. The table highlights how the Farey tree is embedded in butterfly graph as $\phi$ recursions involve Farey sum of four fractions $( \frac{p_R}{q_R}, \frac{p_L}{q_L}, \frac{p_+} {q_+}, \frac{p_-}{q_-}) $ where $\frac{p_{\pm}}{q_{\pm}} = \frac{p_L \pm p_R} { q_L \pm q_R}$. Note that both $\frac{p_+}{q_+}$ and $\frac{p_-}{q_-}$  are neighbors  with the pairs $(\frac{p_L}{q_L},\frac{p_R}{q_R})$, forming friendly triplets.\\
 
 \begin{table}
\begin{tabular}{| c | c |  c  | c | c | }
\hline
Type & $\phi$-recursions  & $\sigma $- recursions &   Example-I &   Example-II \\ 
\hline
\\
$U_L$  & $ \left( \begin{array}{c} \frac{p_R}{q_R}(l+1)  =  \frac{p_{+}}{q_{+}}(l) \\  \frac{p_L}{q_L}(l+1) = \frac{p_L}{q_L} \bigoplus \frac{p_+}{q_+} (l) \end{array}\right) \,\ $ & $ \left( \begin{array}{c}  \sigma_{+}(l+1) \\   \sigma_-(l+1)  \\ \Delta \sigma(l+1) \end{array}\right)   =  \left( \begin{array}{c} \sigma_+(l) + q_L(l)  \\  \sigma_-(l) + q_+(l) 
\\  \Delta \sigma(l)-q_R(l)  \end{array}\right) $   &  $U_L \left( \begin{array}{c} 1 \\ 1 \\  0  \\ \end{array}\right) =\left( \begin{array}{c} 2 \\ 3 \\  -1  \\ \end{array}\right) $ &
$U_L \left( \begin{array}{c} 2 \\ 3 \\  -1  \\ \end{array}\right) =\left( \begin{array}{c} 5 \\ 8 \\  -3  \\ \end{array}\right) $
\\ \hline
\\
$U_R$  & $ \left( \begin{array}{c} \frac{p_L}{q_L}(l+1)  =  \frac{p_{+}}{q_{+}}(l) \\  \frac{p_R}{q_R}(l+1) = \frac{p_R}{q_R}  \bigoplus \frac{p_+}{q_+}(l) \end{array}\right) \,\ $ & $ \left( \begin{array}{c}  \sigma_{+}(l+1) \\   \sigma_-(l+1)  \\ \Delta \sigma(l+1) \end{array}\right)   =  \left( \begin{array}{c} \sigma_+(l) + q_+(l)  \\  \sigma_-(l) + q_R(l) 
\\  \Delta \sigma(l)+q_L(l)  \end{array}\right) $   &  $U_R \left( \begin{array}{c} 1 \\ 1 \\  0  \\ \end{array}\right) =\left( \begin{array}{c} 3 \\ 2 \\  1  \\ \end{array}\right) $ &
$U_R \left( \begin{array}{c} 2 \\ 3 \\  -1  \\ \end{array}\right) =\left( \begin{array}{c} 7 \\ 5 \\  2  \\ \end{array}\right) $ \\ \hline
\\
$D_L$  & $ \left( \begin{array}{c} \frac{p_R}{q_R}(l+1)  =  \frac{p_+}{q_+}(l) \\  \frac{p_L}{q_L}(l+1) = \frac{p_L}{q_L} \bigoplus \frac{p_+}{q_+} (l) \end{array}\right) \,\ $ & $ \left( \begin{array}{c}  \sigma_{+}(l+1) \\   \sigma_-(l+1)  \\ \Delta \sigma(l+1) \end{array}\right)   =  \left( \begin{array}{c} \sigma_+(l) + q_+(l)  \\  \sigma_-(l) + q_L(l) 
\\  \Delta \sigma(l)+q_R(l)  \end{array}\right) $   &  $D_L \left( \begin{array}{c} 1 \\ 1 \\  0  \\ \end{array}\right) =\left( \begin{array}{c} 2 \\ 3 \\  1  \\ \end{array}\right) $ &
$D_L \left( \begin{array}{c} 2 \\ 3 \\  -1  \\ \end{array}\right) =\left( \begin{array}{c} 5 \\ 8 \\  1  \\ \end{array}\right) $ \\ \hline
\\
$D_R$  & $ \left( \begin{array}{c} \frac{p_L}{q_L}(l+1)  =  \frac{p_+}{q_+}(l) \\  \frac{p_R}{q_R}(l+1) = \frac{p_R}{q_R}  \bigoplus \frac{p_+}{q_+}(l) \end{array}\right) \,\ $& $ \left( \begin{array}{c}  \sigma_{+}(l+1) \\   \sigma_-(l+1)  \\ \Delta \sigma(l+1) \end{array}\right)   =  \left( \begin{array}{c} \sigma_+(l) + q_R(l)  \\  \sigma_-(l) + q_+(l) 
\\  \Delta \sigma(l)-q_L(l)  \end{array}\right) $   &  $D_R \left( \begin{array}{c} 1 \\ 1 \\  0  \\ \end{array}\right) =\left( \begin{array}{c} 3 \\ 2 \\  -1  \\ \end{array}\right) $ &
$D_R \left( \begin{array}{c} 2 \\ 3 \\  -1  \\ \end{array}\right) =\left( \begin{array}{c} 7 \\ 5 \\  -4  \\ \end{array}\right) $ \\ \hline
\\
$C_L$  & $ \left( \begin{array}{c} \frac{p_L}{q_L}(l+1)  =  \frac{p_L}{q_L}(l) \\  \frac{p_R}{q_R}(l+1) = \frac{p_L}{q_L} \bigoplus \frac{p_+}{q_+}(l) \end{array}\right) \,\ $ & $ \left( \begin{array}{c}  \sigma_{+}(l+1) \\   \sigma_-(l+1)  \\ \Delta \sigma(l+1) \end{array}\right)   =  \left( \begin{array}{c} \sigma_+(l) + q_L(l)  \\  \sigma_-(l) + q_L(l) 
\\  \Delta \sigma(l)  \end{array}\right) $   &  $C_L \left( \begin{array}{c} 1 \\ 1 \\  0  \\ \end{array}\right) =\left( \begin{array}{c} 3 \\ 1  \\  0  \\ \end{array}\right) $ &
$C_L \left( \begin{array}{c} 2 \\ 3 \\  -1  \\ \end{array}\right) =\left( \begin{array}{c} 8 \\ 3 \\  -1  \\ \end{array}\right) $ \\ \hline
\\
$C_R$  & $ \left( \begin{array}{c} \frac{p_R}{q_R}(l+1)  =  \frac{p_R}{q_R}(l) \\  \frac{p_L}{q_L}(l+1) = \frac{p_R}{q_R} \bigoplus \frac{p_+}{q_+}(l) \end{array}\right) \,\ $  & $ \left( \begin{array}{c}  \sigma_{+}(l+1) \\   \sigma_-(l+1)  \\ \Delta \sigma(l+1) \end{array}\right)   =  \left( \begin{array}{c} \sigma_+(l) + q_R(l)  \\  \sigma_-(l) + q_R(l) 
\\  \Delta \sigma(l)  \end{array}\right) $   &  $C_R \left( \begin{array}{c} 1 \\ 1 \\  0  \\ \end{array}\right) =\left( \begin{array}{c} 1 \\ 3 \\  0  \\ \end{array}\right) $ &
$C_R \left( \begin{array}{c} 2 \\ 3 \\  -1  \\ \end{array}\right) =\left( \begin{array}{c} 2 \\ 7 \\  -1  \\ \end{array}\right) $ \\  \hline
\\
$C_{cL}$  & $ \left( \begin{array}{c} \frac{p_R}{q_R}(l+1)  =  \frac{p_L}{q_L}(l) \\  \frac{p_L}{q_L}(l+1) =  \frac{p_L}{q_L} \bigodot \frac{p_-}{q_-}(l)\end{array}\right) \,\ $& $ \left( \begin{array}{c}  \sigma_{+}(l+1) \\   \sigma_-(l+1)  \\ \Delta \sigma(l+1) \end{array}\right)   =  \left( \begin{array}{c} \sigma_+(l) + q_-(l)  \\  \sigma_-(l) +q_-(l) 
\\  \Delta \sigma(l)  \end{array}\right) $   &  $C_{cL} \left( \begin{array}{c} 3 \\ 1 \\  0  \\ \end{array}\right) =\left( \begin{array}{c} 5 \\ 3 \\  0  \\ \end{array}\right) $ &
$C_{cL} \left( \begin{array}{c} 2 \\ 7 \\  0  \\ \end{array}\right) =\left( \begin{array}{c} 7 \\ 12 \\  -1  \\ \end{array}\right) $ \\ \hline
\\
$C_{cR}$  & $ \left( \begin{array}{c} \frac{p_L}{q_L}(l+1)  =  \frac{p_R}{q_R}(l) \\  \frac{p_R}{q_R}(l+1) = \frac{p_R}{q_R}  \bigodot \frac{p_-} {q_-}(l) \end{array}\right) \,\ $ & $ \left( \begin{array}{c}  \sigma_{+}(l+1) \\   \sigma_-(l+1)  \\ \Delta \sigma(l+1) \end{array}\right)   =  \left( \begin{array}{c} \sigma_+(l)+ q_-(l)  \\  \sigma_-(l) + q_- (l)
\\  \Delta \sigma(l)  \end{array}\right) $   &  $C_{cR} \left( \begin{array}{c} 3 \\ 1 \\  0  \\ \end{array}\right) =\left( \begin{array}{c} 3 \\ 5 \\  0  \\ \end{array}\right) $ &
$C_{cR} \left( \begin{array}{c} 8 \\ 3 \\  0  \\ \end{array}\right) =\left( \begin{array}{c} 13 \\ 8 \\  -1  \\ \end{array}\right) $ \\  \hline
\end{tabular}
\caption{ Table lists eight generators that act on the column matrix made up three integers  $(q_R(l), q_L (l) , \Delta \sigma (l) )$ to determine the sub-structure of a butterfly. This is illustrated with two examples in the last two columns. The $\phi$ recursions reveal how the Farey tree is embedded in the octonary tree. Note that the symbol $\bigodot$ is analog to $\bigoplus$ where ``plus" is replaced by ``minus" ( $\frac{p_L}{q_L} \bigodot \frac{p_R}{q_R} = \frac{p_-}{q_-}$ ). }
\label{T1}
\end{table}

All matrices representing the eight generators have unit determinant. Furthermore, the $C$ matrices have eigenvalues unity while the $U$ and $D$ matrices have two of the eigenvalues equal to $\zeta=\frac{3 \pm \sqrt{5}}{2}$.  It describes parity non-conservation hierarchy as reflected in the denominators of the as rational approximants of $\zeta$ which can be both even or odd integers. We note that all the self-similar hierarchies correspond to periodic paths in the tree. Described by the product of various combinations of eight generators, their eigenvalues give the self-similar scaling exponents of various butterfly hierarchies. It turns out
that all scaling exponents are
 numbers of the form $\Big (\frac{n+1}{2}+ \sqrt{(\frac{n+1}{2})^2-1} \,\ \Big )$. In continued fraction representation, they are given by  $[(n+1); \overline{1,n}]$. Here $ n=1, 2, 3...$.  This intriguing result can be verified by forming various products of the eight matrices. We recall that
this set of scalings was obtained  explicitly in our previous studies, by an elegant  renormalization equation\cite{SW,SAT21}.

\subsection{ The Butterfly Tree }

  \begin{figure}[htbp] 
 \includegraphics[width = .35 \linewidth,height=.6 \linewidth]{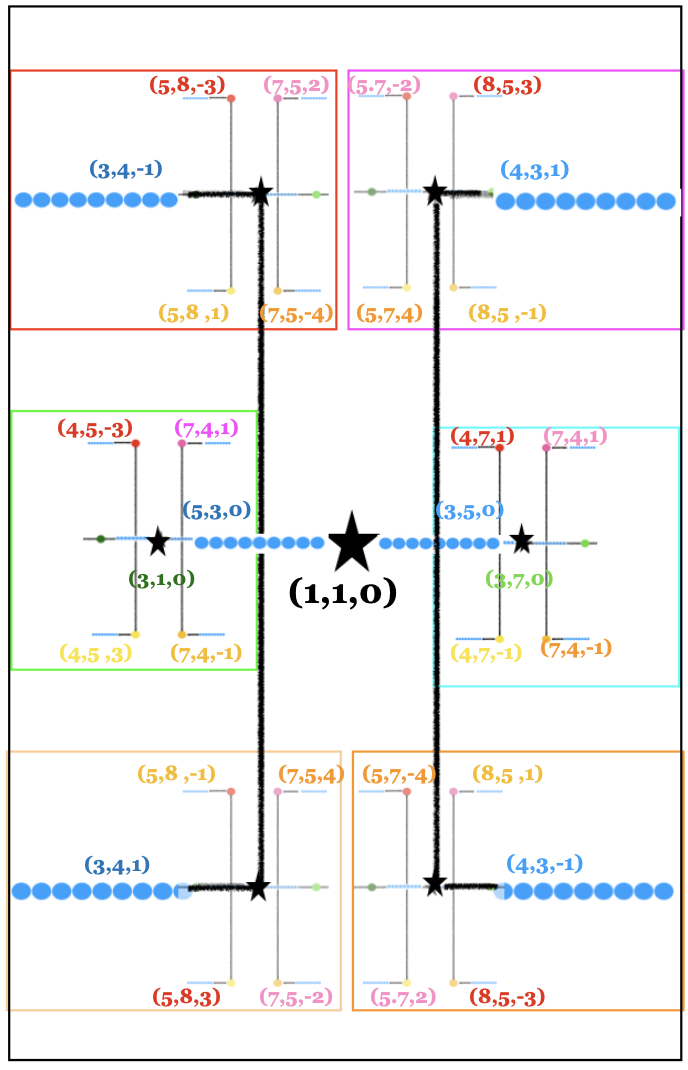} 
\caption{ Eightfold butterfly tree where parent butterflies are shown with a star. Figure shows three generations starting with main butterfly. In the second generation the babies themselves become parent butterfly.  Sextuplets are shown in six boxes ( color codes ), each with a tail shown with a row of blue dots.
Butterflies in the chain also follow the same hierarchical eightfold pattern which is not shown in figure. Tree is labeled with the triplet $(q_R, q_L, \Delta \sigma)$
which is color coded.  }
\label{Btree}
\end{figure}

Figure (\ref{Btree}) shows the octonary butterfly tree, color coded, labeled with three integers $(q_R, q_L, \Delta \sigma)$. The outer big box houses the main butterfly ( labeled in black with $(1,1,0)$  ) with
six babies, each with a tail - an infinite chain of butterflies. 
The sub-structure of each of the six butterflies are shown in  smaller boxes. Blue integers show the labelling of the first member of the chain that is attached to every butterfly.\\

To summarize the iterative process, one starts with a butterfly and use $( C_L, C_R, U_L, U_R, D_L, D_R)$ generators to produce the sextuplets. One then uses $C_{cL}$ for $(C_L, U_L, D_L)$ generated butterflies and $C_{cR}$ for $(C_R, U_R, D_R )$ produced butterflies to generate their tails.
The process is repeated ad infinitum for each of the six butterflies as well as for the chains  consisting of infinity of of butterflies that form the tails.  The resulting
Butterfly tree  is isomorphic to the butterfly fractal. Labelling of the tree with three integers provides a unique label for every butterfly in the fractal graph as no two butterflies share  the same label. 
The topological parameter $\Delta \sigma$ encoding butterfly asymmetry emerges as a key parameter in the butterfly graph, whose values exhibit a very orderly and predictable pattern in the tree. This is in sharp contrast to $(\sigma_+, \sigma_-)$ labels exhibiting a rather erratic pattern. We summarize some of the key features of the tree.

\begin{itemize}
\item C-cell of baby butterflies have same $\Delta \sigma$ as their parent butterfly.
\item For E-cell butterflies, $\Delta \sigma_L + \Delta \sigma_R = q_c$. Here $ \Delta \sigma_L$ and $ \Delta \sigma_R$ denote the asymmetry parameter for the left and the right edge butterflies and   $q_c$ is the denominator of the rational flux for the center of the parent butterfly.
\item Chains are considered as part of the C-cells as all butterflies in the chain have same parity and also same $\Delta \sigma$.
\item Chain of butterflies attached to every butterfly are like a tail of the butterfly with a fixed $\Delta \sigma$ characterizing the entire chain which is the $\Delta \sigma$ of the parent butterfly.
In other words, $\Delta \sigma$ is the most important {\it gene} of a butterfly and the entire chain attached to it inherits that gene.
\end{itemize}
\section{ Butterfly  Fractal and Integer Apollonian Gasket }
Relationship between the butterfly fractal and Integer Apollonian gasket- an abstract geometrical fractal describing packing of mutually tangent circles in a plane, has been subject of earlier studies\cite{book,SAT21}.
The four generators of the Apollonian group that describe the packing of mutually tangent circles are:
\begin{equation}
S_1  =  \left( \begin{array}{cccc} -1 & 2 & 2  & 2 \\    0 & 1 & 0 & 0  \\  0 & 0 & 1 & 0 \\ 0 & 0 & 0 & 1  \\ \end{array}\right),\,\ 
S_2  =  \left( \begin{array}{cccc} 1 & 0 & 0  & 0 \\    2 & -1 & 2 & 2  \\  0 & 0 & 1 & 0 \\ 0 & 0 & 0 & 1  \\ \end{array}\right),\,\
S_3  =  \left( \begin{array}{cccc} 1 & 0 & 0  & 0 \\    0 & 1 & 0 & 0  \\  2 & 2 & -1 & 2 \\ 0 & 0 & 0 & 1  \\ \end{array}\right),\,\
S_4  =  \left( \begin{array}{cccc} 1 & 0 & 0  & 0 \\    0 & 1 & 0 & 0  \\  0 & 0 & 1 & 0 \\ 2 & 2 & 2 & -1  \\ \end{array}\right)
\end{equation}

As pointed out in previous study\cite{SAT21}, one needs
super Apollonian group that  also includes the adjoint of $S_i$ to describe butterfly fractal. They are needed to accommodate both the upper and the lower edges of a butterfly. 
We now sharpen the comparison between the butterfly fractal and the Apollonian gasket by relating the generators that describe them.\\

In view of the correspondence between Ford Apollonian (subset of Apollonian where one of the circle is the $x$-axis ) and the Pythagorean triplets\cite{SAT21}, it is not surprising that the three operators $(H_1, H_2, H_3)$ form a subset of the four operators $(S_1, S_2, S_3, S_4)$ that
underlie Apollonian packing.
Straight forward analysis reveal the following correspondence between the $H_i$ and the $S_i$:
\begin{eqnarray*}
H_1  & \leftrightarrow &   S_1 S_2, \,\,\,\,\,\,\,\,\,\ H_2 \leftrightarrow S_1 S_3 ,\,\,\,\,\,\,\,\,\,\,\  H_3 \leftrightarrow  S_1 S_3 S_1,\,\,\,\,\,\,\,\,\,\,\  \
 U_L \leftrightarrow  S_4 S_2 ,\,\,\,\,\,\,\,\,\,\,\  \ U_R \leftrightarrow  S_4 S_3  \\
\label{HS}
\end{eqnarray*}

The four adjoint operators provide the correspondence with the $(D_L, D_R)$ operators.  We note that  in contrast to  $H_i$ operators, $S_i$ operators allow parity non-conserving hierarchies.\\

Both the Butterfly group and the super Apollonian group have eight generators. However, inspite of all its elegance, the super Apollonian group does not provide one to one mapping between the Pythagorean quadruplets and the butterfly fractal. This is due to
edge-edge siblings as was pointed out in our earlier study\cite{SAT21}. We believe that the problem can be remedied by replacing circles by spheres whose great circles lie in the 2D plane
and the degenerate configurations can be described by out of plane great circles and this is a work in progress.\\

\section{Summary and Conclusion}

This is the first paper that reveals the tree structure of the butterfly fractal. Labeled with integers, it is one to one mapping with the hierarchical structure of the butterfly graph that provides
an iterative scheme to build the butterfly fractal.
 Novel feature in this formulation is to include the renormalization equations for  $\sigma_{\pm}$,  the topological quantum numbers for the gaps of the spectrum. These quanta are experimentally accessible quantities and are key to describe topological aspects of a butterfly, something that was missing in the 
earlier studies\cite{SW,SAT21} . We want to emphasize that it is the gap Chern number that are indispensable in building the butterfly graph  and not the Band Chern numbers.\\
 
  Identifying butterfly with a tail as the building block of the butterfly fractal leads to the theoretical framework of eightfold way to build this quantum fractal.
  The Butterfly group  with eight generators represented by matrices with integer coefficients is a generalization of the modular group of $ 2 \times 2$ matrices with unit determinant and integer coefficients denoted as $SL(2,Z)$ group.  Underlying this generalization is the inclusion of $\Delta \sigma_{\pm}$ 
that quantifies the distortions in the sub-images in the butterfly graph.  Central cells of the butterflies preserve this entity at all scales as  $\Delta \sigma$ remains invariant under iteration.\\


The hierarchical sets of C-cell butterflies provide a physical realization of the Pythagorean tree. So these butterfly hierarchies are also the hierarchies of the Pythagorean trees 
 each labeled with an integer $\Delta \sigma$, embedded in the butterfly graph. The geometrical interpretation of the ternary nature of the Pythagorean tree have been described  earlier\cite{PT3}.
 Butterfly graph complements this geometrical picture as  Nature has found a beautiful way to display the entire set of primitive set Pythagorean triplets in every sub-images. It is an interesting open question whether eightfold tree has an abstract representation in mathematics of abstract hierarchical sets.\\
 
 
The repertoire of this paper clearly shows the importance of number theory in describing this quantum fractal. The Diophantine way of labeling butterflies,
as it emerges in our theoretical framework proves the fact that underlying butterfly fractal is an integer labelling of the fractal. In other words, butterfly fractal is made up of integers which determine the hierarchical pattern. Quantum mechanics determine the widths of the energy bands, the only feature of the butterfly graph that falls outside the number theoretical framework.\\

\end{document}